\documentclass[10pt]{amsart}
\usepackage{amsmath,amssymb,amsfonts,amsthm,amstext,verbatim}
\usepackage{enumerate,color,graphicx,stmaryrd,psfrag} 
\usepackage{url}
\usepackage{mdwlist}   
\usepackage{mathtools} 
\usepackage{wasysym}   

\RequirePackage[colorlinks,citecolor=blue,urlcolor=blue]{hyperref}
\newcommand\oo{\infty}

\newcommand\g{\gamma}

\newcommand\pc{p_{\text{\rm c}}}

\newcommand\resp{respectively}

\newcommand{\ZZ}{\mathbb{Z}}     

\def\mik{1}
\newcommand\cpsfrag[2]{\ifnum\mik=1\psfrag{#1}{#2}\fi}

\numberwithin{equation}{section}
\numberwithin{thm}{section}
\numberwithin{figure}{section}

\newcounter{mycount}

\renewcommand{\subparagraph}[1]{}

\newcommand\sL{{\mathcal L}}

\begin{document}
\title{Percolation since Saint-Flour}

\author[Grimmett]{Geoffrey R.\ Grimmett}
\address{Statistical Laboratory, Centre for
Mathematical Sciences, Cambridge University, Wilberforce Road,
Cambridge CB3 0WB, UK} 
\email{\href{mailto:g.r.grimmett@statslab.cam.ac.uk}{\nolinkurl{g.r.grimmett@statslab.cam.ac.uk}}}
\urladdr{\url{http://www.statslab.cam.ac.uk/~grg/}}

\author[Kesten]{Harry Kesten}
\address{Department of Mathematics, Malott Hall, Cornell University,
Ithaca, NY 14853-4201, USA}
\email{\href{mailto:kesten@math.cornell.edu}{\nolinkurl{kesten@math.cornell.edu}}}
\urladdr{\url{http://www.math.cornell.edu/People/Faculty/kesten.html}}



\maketitle

\section{Introduction}

There has been a great deal of interest and activity 
in percolation theory since the
two Saint-Flour courses, \cite{G-stf,K-stf}, of 1984 and 1996
reprinted in this new edition. 
We present here a summary of progress since the first publications of our lecture notes.
 
The second edition of  \emph{Percolation},
\cite{Grimmett_Percolation}, was published in 1999 as a fairly
complete, contemporary account of `classical' percolation
(but with only limited coverage of first-passage percolation).
Books published since in the general area include \cite{Grimmett_RCM,Grimmett_Graphs} 
and, for two-dimensional
theory, \cite{BolRio}.   There has been in addition a third Saint-Flour
course on percolation, on the Wulff construction by Rapha\"el Cerf in 2004, 
\cite{Cerf-Wu}. The 2004 review by Howard of first-passage percolation, \cite{How04}, 
is useful. 

The principal targets of this short introduction are to indicate areas of significance that have received
attention in recent years, and to include a useful bibliography.  
It is not intended as a comprehensive survey of the vast amount
of material related to percolation and published since 1984/1996.
We restrict our attention mostly to problems associated with 
finite-dimensional lattices, and we take quite a broad view of 
what constitutes a percolation or first-passage percolation problem.
We offer apologies to the many authors whose work does not
receive explicit mention here, and we trust sufficient pointers
are present to guide the reader through recent developments in percolation
theory.  

\section{Aspects of first-passage percolation}

\subsection{Limit theorems}

The main `ergodic theorems' for first-passage percolation
were established in their original forms by Hammersley and Welsh in 1965, and 
in more refined forms by others later; see \cite{K-stf} for the history prior to 1984.
It is an important open problem to identify and prove the corresponding distributional limits.
Only a limited amount is known beyond the Gaussian limit theorem 
of Kesten and Zhang
for a critical model in two dimensions, \cite{KesZ}. 
They were able to exploit the special structure of dual cycles
in critical percolation.

One  expects a Gaussian limit for
passage-times along a long tube, but how wide can such a tube be?
In recent work \cite{ChatD}, Chatterjee and Dey have shown a Gaussian limit
for $d$-dimensional tubes with width growing as a power of length
(see also \cite{Ahl}). 
There is evidence that the largest acceptable power is connected to the 
so-called fluctuation exponent
(see Section \ref{sec:kpz}).

The weak limit is not expected to be Gaussian for 
unrestricted first-passage problems (see the corresponding discussion of
Section \ref{sec:lpp} for last-passage percolation).
Estimates of variance are first 
steps towards identification of distributional limits. Following earlier
results by others, Benjamini, Kalai, and Schramm, \cite{BenKS},
showed  that the variance grows no faster than $n/\log n$
in a special  first-passage process.
Their methods and conclusions have been extended by Bena\"{\i}m and Rossignol, \cite{BenR}.

Large-deviations for passage-times have been discussed since \cite{K-stf} by Chow and Zhang, \cite{ChowY}, and  Cranston, Gauthier, and Mountford, \cite{CGM}. 
Corresponding results for maximum flows are summarized in Section \ref{sec:mf}.

LaGatta and Wehr have initiated a theory of  `Riemannian' first-passage percolation.
The shape theorem of \cite{LGW10} has been followed by a more detailed study of
geodesics, \cite{LGW12a,LGW12b}.

\subsection{Fluctuation and wandering exponents}\label{sec:kpz}

Some aspects of the geometry of first-passage percolation are not yet well understood.
For example, how great are the fluctuations in fastest paths, when
do asymptotic shapes have positive curvature, how many infinite geodesics exist?   
These intertwined questions have attracted a great deal of attention
in recent years. We sketch certain
recent advances concerning a fluctuation theory for fastest paths.

For two vertices $x$, $y$, write $T(x,y)$ for the passage-time from $x$ to $y$,
and $D(x,y)$ for the maximum deviation between the 
fastest path from $x$ to $y$ and
the straight line-segment from $x$ to $y$. 
For clarity, we assume that the passage-times 
have a continuous distribution, so that
the fastest path between two given vertices is almost surely unique.

Roughly speaking, the \emph{fluctuation exponent} is the number
$\chi$ such that $T(x,y)$ deviates from its mean by order $|x-y|^\chi$. The \emph{wandering
exponent} is the number $\xi$ such that $D(x,y)$ is typically of order $|x-y|^\xi$.
The universal relation 
$\chi=2\xi -1$ is believed to hold, 
and furthermore $\chi=\frac13$, $\xi=\frac23$
in two dimensions. This relation is sometimes referred to
as the KPZ relation, after the authors of \cite{kpz} (see also \cite{KrugS}).

Several possible definitions of the exponents have been discussed in the literature.
Following earlier progress by others including Newman and Piza, \cite{NewPiz},
Chatterjee has succeeded in proving a version of the KPZ relation,
under suitable definitions and subject to an assumption on the
passage-time distribution, see \cite{Chat11}.
His conclusions have been refined and improved by 
Auffinger and Damron, \cite{AuffD11}.

\subsection{Geodesics}

A \emph{geodesic} is a semi- or doubly-infinite path $\pi$ such that, for all
vertices $u,v \in \pi$, the section of $\pi$ from $u$ to $v$ is shortest in
the metric induced by the passage-times. 
Motivated by work with  Stein on spin glasses,
Newman asked in \cite{Newman95} for bounds on the number
of geodesics that can exist in a $d$-dimensional first-passage process. 
He has conjectured that, subject to natural conditions,
no doubly-infinite geodesics (`bi-geodesics') exist, and there are infinitely many
semi-infinite geodesics (`uni-geodesics').  Existence of
bi-geodesics on $\ZZ^2$ is equivalent to existence of non-constant ground states
for the related disordered ferromagnet, \cite{New97}.

The existence (or not) of geodesics is bound up with the curvature of the boundary of the limiting shape
of the growth process, and also with the question of coexistence of
populations in a competition model. Results to date are incomplete. 
Interested readers are directed to the recent paper of Auffinger and Damron,
\cite{AuffD11b}, and to \cite{BS10,GarM, Hoff08,Lall,LicN,We97,We98}
and the references therein.

The geometry of geodesics in planar \emph{Riemannian} first-passage percolation is the subject
of a recent project of LaGatta and Wehr, \cite{LGW12a,LGW12b}.

\subsection{The maximum flow problem}\label{sec:mf}

In two dimensions, the dual of first-passage percolation is the problem
of finding the maximum flow in a certain (dual) random network. 
The maximum-flow problem is well posed in higher dimensions, and is dual to
the problem of finding the minimum size of certain `cut-surfaces'.
The last problem is much harder than the more usual first-passage problem,
since  tricky geometrical issues arise in stitching
surfaces together. 

Partial progress was made by Kesten, \cite{Kes87}, as reported
in \cite{K-stf}. A number of authors have worked on the problem since, 
with substantial results by Zhang, \cite{Zh,Z07}, 
Rossignol and Th\'eret, \cite{RossT10a,RossT10b},
and Th\'eret, \cite{Th07,Th08}.  Remarkable progress has been made over
a period of years, culminating in a recent series of papers
by Cerf and Th\'eret, \cite{Cerf11, CerfT-LD1, CerfT-LD2}.
They consider the general situation of maximal flows in a scaling limit of $n^{-1}\ZZ^d$ as $n \to\oo$.
In the special case of a cube in $d$ dimensions, these results yield that
the maximum flow between opposite faces of
the $n$-cube is asymptotic to $\phi n^{d-1}$
for some $\phi = \phi(F)$, where $F$ is the distribution of a typical edge-capacity.
In addition, the authors of these papers have proved lower and upper large-deviation
theorems  (of surface order $n^{d-1}$, and volume order $n^d$, \resp), and 
also that $\phi>0$ if and only if the atom $F(0)$ is sufficiently small that
the edges with weight $0$ do not percolate.

\subsection{Last-passage percolation}\label{sec:lpp}

Johansson, \cite{Joh00,Joh08}, showed the presence of  the Tracy--Widom distribution
within a certain percolation-type problem. While this is a \emph{last-} rather
than a \emph{first}-passage problem, it is
included here both for its intrinsic interest, and for its 
potential implications
for weak limits of \emph{first}-passage times. To each edge 
$e$ of the square lattice $\ZZ^2$ we allocate
a random variable $T_e$ with the geometric 
distribution, $P(T=k) = (1-q)q^k$, $k \ge 0$.
Let $G_{m,n}$ be the supremum of the passage-times of
paths from $(1,1)$ to $(m,n)$ each of whose steps are either upwards or rightwards, and 
let $\gamma > 1$. It is shown in \cite{Joh00} (see also \cite{Joh02}) 
that there exist functions $\mu(\g,q)$,
$\sigma(\g,q)$, with known closed forms,  
such that $(G_{\g n,n} - \mu n)/(\sigma n^{1/3})$ is
distributed asymptotically as the Tracy--Widom distribution. 
The exponent $\frac13$ occurs here just as in the problem of the longest increasing 
subsequence, and this is the conjectured value of the 
fluctuation exponent $\chi$ in two dimensions.

Johansson's work develops several beautiful connections:
to so-called determinantal processes,  \cite{Joh06},
to the Kardar--Parisi--Zhang (KPZ) universality class, \cite{Cor11,Joh06,kpz,KrugS}, and  to
the `totally asymmetric simple exclusion process' TASEP, \cite{BA-Cor}.

Graham, \cite{Gr10},  has explained how to use the theory of concentration
to obtain the sublinearity of variance for last-passage percolation.
We mention also the recent papers of: Hambly and Martin, \cite{HamM},
who consider edge-weights with heavy tails;
Lin and Sepp\"al\"ainen, \cite{LSep}, concerning limiting shape;
and Kesten and Sidoravicius, \cite{KS07}, concerning directed last-passage percolation.

\subsection{Greedy lattice animals}

To the vertices  of a lattice we allocate random real-valued scores.
A \emph{greedy lattice animal} of size $n$ is a connected subgraph
of $n$ vertices including the origin such that its aggregate score is 
maximal across this class. This variant of first-passage percolation
was initiated by Cox, Griffin, Gandolfi and Kesten, \cite{CGGK,GKes94}.
A connection with Euclidean first-passage percolation was explored by
Howard and Newman, \cite{HNew99}.
The geometry of greedy animals has been investigated by Hammond, \cite{Ham06},
where further references may be found.

\section{Percolation and disordered systems}

\subsection{Conformality in two dimensions}\label{sec:c2d}
Possibly the most prominent recent advance in probability theory 
is Schramm's construction of the family
of random curves now termed Schramm--L\"owner curves (SLE);
see \cite{Sch00,Sch06}.
The theory of SLE has
revolutionized the mathematics of planar phase transitions, and is a pillar of
the bridge between probability and conformal field theory.

In applying SLE to percolation, a
spatial limit theorem is required in order to transform
the critical process on a discrete grid into a process on a continuum.
This limit theorem has been proved in essentially only
one case, namely site percolation on the triangular lattice.
A key ingredient for percolation was the 
so-called Cardy formula of \cite{Cardy}, which is
based on the arguments of conformal field theory, and asserts
a formula for the probabilities of open crossings of large domains in critical models.
Cardy's formula was proved in 2001 by Smirnov, \cite{Smirnov,SmirnovII},
for the site model on the triangular lattice. The
last paper included an outline proof of what has been called `full conformal invariance',
and this has been fulfilled and amplified since
by Camia and Newman, \cite{CamNew,CamNew2}.

It is a significant open problem to prove 
conformal invariance
for an extensive family of critical two-dimensional percolation processes.  

The theory of SLE has been described in the Saint-Flour notes
of Wendelin Werner, \cite{W-stf}, as well as in
Greg Lawler's book \cite{Lawler}, the recent surveys 
\cite{Beffara2010,BD-rev,law11,sun},
and the references therein. Proofs and elaborations of Cardy's formula have appeared in 
\cite{BolRio,Cam-rev,Grimmett_Graphs,ShW,WW_park_city}.

Camia, Fontes, and Newman have asked about scaling limits for \emph{near}-critical percolation,
\cite{CFN1,CFN2}. A series of papers on this topic has been 
promised by Garban, Pete, and Schramm, of which the first is \cite{GPS12}.

In contrast to the case of two dimensions, only little
has been achieved  since \cite{G-stf,Grimmett_Percolation}
for critical percolation in (slightly) higher dimensions. 
For example, it remains a tantalizing problem to show
the continuity of the percolation probability \emph{at} the
critical point. Similarly, there has been essentially no progress
on the percolation phase transition in $d$ dimensions with $3 \le d \le 18$.
When $d \ge 19$, we are in the domain of the lace expansion, 
see \cite[Sect.\ 10.3]{Grimmett_Percolation}.

\subsection{Criticality and universality in two dimensions}

A core technique for the study of near-critical percolation in two dimensions
is the theory of box-crossings: when do they exist, and how may they be combined?
The methods of Russo and Seymour--Welsh (RSW) have proved very useful
in a variety of settings, including the proof of Cardy's formula
and conformality. More
recently, they have been used to prove aspects of universality for bond percolation on
isoradial graphs. This large family of models includes but is much broader than
the more familiar inhomogeneous square, triangular and hexagonal models.
See \cite{GM1,GM2,GM3} and the overview \cite[Sect.\ 3]{G-three}.

RSW theory has become streamlined by the recent introduction of 
the theory of influence and sharp threshold, reviewed
in \cite[Sect.\ 4.5]{Grimmett_Graphs}. Examples,
in addition to the universality results above, include the
proofs by Bollob\'as and Riordan that the critical probability
of  site percolation on the Voronoi graph of a Poisson process is $\frac12$,
\cite{BR06b, BR08}, and their exploration of the relationship
between self-duality and criticality for certain polygon models
in two dimensions, \cite{BR}.

Mention is made here of the recent proof by Schramm and Smirnov, \cite{SS11},
that any scaling limit of a planar percolation is a so-called `black noise'.

The \emph{connective constant} of a lattice $\sL$ is the number $\kappa$
such that the number of $n$-step self-avoiding walks on $\sL$ grows in the
manner of $\kappa^n$. Identifying the value of $\kappa$
is akin to calculating the value of a percolation critical point. 
Nienhuis's long-standing conjecture that $\kappa = \sqrt{2+\sqrt 2}$
for the hexagonal lattice 
has been proved in notable work of  Duminil-Copin and Smirnov, \cite{Dum-S}.
It is an important open problem to show that the weak limit  
of a random self-avoiding walk in two dimensions is the
Schramm--L\"owner curve SLE$_{8/3}$
(see \cite{BDGS}).

\subsection{Dynamical percolation and noise sensitivity}\label{sec:dyn}

Suppose a critical percolation model is allowed to evolve in time: the state of each bond/site
switches between open and closed after exponential holding times, with rates chosen 
in such a way that the
process remains critical. This `dynamical percolation'
process was introduced in independent work of H\"aggstrom, Peres, and Steif,
\cite{HPS97}, and of Benjamini.  
The main issue has been to determine the Hausdorff dimensions
of the sets of times at which certain forms of aberrant behaviour occur. 
It has been normal to consider the case of two dimensions, and
particularly (but not exclusively)
site percolation on the triangular lattice, since
a certain amount is then known about pivotal sites for long open connections.
Earlier results of Schramm and Steif, \cite{SSt}, for the existence of an infinite
open cluster have been 
complemented by a more detailed analysis by Garban, Pete, and Schramm, \cite{GPS11}.
Dynamical percolation has been reviewed recently by Garban and Steif, 
\cite{GarS}; see also \cite{BGS,Gar11-S,St09}.

In the related area of `noise-sensitivity' of Boolean functions, 
Benjamini, Kalai,
and Schramm  have asked how sensitive is the probability of an event 
to the resampling of the states of bonds/sites, \cite{BKS99}.
Noise sensitivity has been reviewed  by Garban, \cite{Gar11-S}, and
jointly with Steif, \cite{GarS}.

\subsection{Incipient infinite cluster}\label{sec:iic}

As remarked in Section \ref{sec:c2d}, 
it remains an open problem to prove that there is no infinite cluster 
in the \emph{critical} percolation model in a general number of dimensions.
This is known when $d=2$ (by box-crossing arguments) and when $d \ge 19$ (using
the lace expansion). The so-called `incipient infinite cluster' (IIC) is a measure
defined via conditioning at or near the critical point, under which an infinite cluster
exists. Kesten's proof of existence of the IIC in two dimensions has been extended 
and developed by J\'arai,
\cite{jar}. Corresponding results in high dimensions 
have been established by van der Hofstad and J\'arai, \cite{vdHJ},
using the lace expansion, and more recently by Heydenreich, van der Hofstad, and Hulshof,
\cite{HvdHH}, assuming a form of the triangle condition.
The latter results are valid also for sufficiently spread-out models in $7$ and more dimensions.

Much is known about the IIC for oriented percolation, see Section \ref{sec:orient}.

\subsection{Random walk on the infinite cluster}

Supercritical percolation possesses an almost surely unique infinite cluster.
This cluster has many properties in common with the entire lattice, and
these two sets may be said to be `coarsely' equivalent. One test
for coarse equivalence is whether or not random walk on the infinite cluster satisfies
an invariance principle. A number of authors have worked on this problem, with
notable success. No complete bibliography is attempted here. We refer the reader instead to
the survey by Biskup, \cite{Bisk11}, of the more general topic
of random walks on lattices endowed with random conductances.
 
There is an old conjecture of Alexander and Orbach for the rate of decay
of the return probability of random walk on the incipient infinite cluster,
namely that $p_{2n}(x,x)$ is of order $n^{-2/3}$. This has been proved by Kozma and
Nachmias, \cite{KozNac}, for large $d$ ($d \ge 19$ suffices) and for 
the spread-out model with $d > 6$. A related result for oriented percolation
is found in Section \ref{sec:orient}.

We mention also recent work of Fribergh and Hammond, \cite{FriH}.
In a continuation of earlier work by others, they
have studied the velocity of a random walk with drift on the supercritical infinite open 
cluster. They have proved the (arguably counter-intuitive) 
existence of a critical value for the 
drift rate below which the velocity is strictly positive, and above which
it is $0$.

\subsection{Subcritical percolation}

In the subcritical phase, the two-point connectivity function decays exponentially with distance $n$.
That is, the probability of an open connection between two points at distance $n$ 
behaves qualitatively in the manner of $e^{-n/\xi(p)}$ where $\xi(p)$ is the `correlation length'.
There is a multiplicative correction term $n^{-(d-1)/2}$, and the
ensuing decay is termed `Ornstein--Zernike decay' (OZ), \cite{OZ}.
This polynomial correction 
has been proved throughout the subcritical phase by Campanino and Ioffe, \cite{CamI},
in a refinement of earlier work restricted to axial directions.  
In the supercritical phase, one considers instead the \emph{truncated} two-point function,
for which the most recent study is that of Campanino and Gianfelice, \cite{CamG}.
It is an open question to prove OZ decay for the truncated function all the 
way down to the critical point. 

\subsection{Supercritical percolation}

A finite open cluster in two dimensions is surrounded by a dual cycle.
Duality is much more complex in  dimensions exceeding $2$,
since the boundary of a finite cluster is a surface, and general surfaces  have complicated 
geometry. The geometry of random surfaces has been the source of a number of fine problems,
of which we mention three.

Whereas classical percolation concerns the binary relation of
`connectivity', one may instead seek conditions for there to
be an infinite `entangled' set of open edges (in three dimensions). 
Unlike `connectivity', entanglement is a long-range property, and its
study involves novel complications  of a geometrical nature. 
Early papers \cite{GHol00,Hag01e,H00}
by Grimmett, H\"aggstr\"om, and Holroyd 
stimulated recent work of Atapour and Madras, \cite{atapour-madras}, who have 
proved an exponential-tail theorem for the size of the largest entanglement at the origin
in part of the subcritical phase. 
An improved lower bound for the critical 
entanglement probability has been proved by Grimmett and Holroyd, \cite{GHol10}.

The last reference is the second of a sequence of four papers by the same authors,
sometimes joined by Dirr, Dondl, and Scheutzow, devoted to 
questions associated with
so-called `Lipschitz percolation', see \cite{DDGHS, GH10c,GH10b}.
Lipschitz percolation is concerned with the existence and geometry of 
`flat' open  surfaces within
a percolation process.

One may ask similarly for conditions under which there exists an infinite `rigid' set of
open edges. This problem has been studied by Holroyd, 
\cite{H01r}, and H\"aggstr\"om, \cite{Hag01r}.  In \cite{CRV,MRV}, Connelly, Menshikov,
Rybnikov, and Volkov have proved loss of tension in planar networks  subjected
to random perforations.  

By conditioning on the absence of open connections between the upper
and lower hemispheres of  a large sphere, Gielis and Grimmett, \cite{GG02},
have constructed Dobrushin interfaces for percolation in three and more dimensions, 
and more generally for random-cluster models. The associated `roughening transition' is
no better understood for such systems than for the Ising model.

The \emph{chemical distance} between two points is defined to be
the length of the shortest open path between these points. Antal and Pisztora
showed the `law of large numbers' for chemical distance, \cite{AP}.
Their work has been elaborated recently by Garet and Marchand in their studies
of large and moderate deviations, \cite{GM07,GM10}.

\subsection{Oriented percolation}\label{sec:orient}

The existence of the incipient
infinite cluster for oriented percolation 
has been proved in $d > 4+1$ dimensions by
van der Hofstad, den Hollander, and Slade, \cite{vdHHS}.

In  a useful paper concerning random walk on random graphs, Barlow, J\'arai,
Kumagai, and Slade, \cite{BJKS}, have proved a version of the Alexander--Orbach
conjecture in the setting of the incipient infinite cluster
of oriented percolation in $d >7$ dimensions.

The geometrical structure of supercritical oriented percolation has been investigated 
by Grimmett and Hiemer, \cite{GHiemer}, 
and further in two dimensions by Wu and Zhang, 
\cite{WuZ}.

\subsection{Invasion and bootstrap percolation}

`Invasion percolation' is a system that is \lq\lq self-organized" to
focus on the critical point of percolation.  It proceeds
by a progressive filling of so-called ponds. Quite  a lot
is now known about these ponds and their outlets.
The reader is referred to the recent papers \cite{DamS11,DS12} by 
Damron and Sapozhnikov, whose results have been inspired in part by earlier work
of others, \cite{AGHS,Good}, on invasion percolation on trees.

In `bootstrap percolation', each vertex in a box of side-length $n$ is declared
`infected' with probability $p$, and further vertices are infected if they have two
or more infected neighbours.  
The \emph{critical point} $\pc=\pc(n)$ is the value of $p$ above which there is probability
at least $\frac12$ that every vertex becomes infected. 
Increasingly accurate asymptotics
of $\pc(n)$, viewed as a function of $n$ as $n \to\oo$, have been obtained.
Gravner, Holroyd, and Morris, \cite{GHM},  have provided the second-order term in the expansion
of $\pc(n)$ in two dimensions, thereby complementing the earlier result of Holroyd
that $\pc(n) \sim \pi^2/(18\log n)$. The sharpness of
the phase transition for a general class of such models has been established
by Balogh, Bollob\'as, Duminil-Copin, and Morris, \cite{BBDM}.

\subsection{Percolation of words}

Each vertex of a lattice is allocated a random letter from the two-letter alphabet
$\{0,1\}$. A \emph{word} is a semi-infinite sequence of letters. Benjamini and Kesten,
\cite{BenK}, initiated the study of the set of words seen in a two-dimensional lattice. 
Their work has been continued by Kesten, Sidoravicius, and Zhang, \cite{KSZ}.

Grimmett and Holroyd, \cite{GH10b}, have presented related results
for Lipschitz embeddings of $d$-dimensional words in $d'$-dimensional space. 
The case $d=d'=1$ was posed by Grimmett, Liggett, and Richthammer, \cite{GLR},
and has been largely solved in recent work of Basu and Sly, and of G\'acs, \cite{BSly,gacs}. 

\subsection{Percolation on non-amenable graphs}

In the last 20 years, there has been a good deal of interest in random processes
on graphs whose surface/volume ratio does not
tend to zero. Grimmett and Newman, \cite{GrimNew},
proved the existence of three phases for percolation on the direct product 
of a tree and a lattice. This provoked Benjamini and Schramm, \cite{BenS96},
to propose in 1996 a systematic study of
percolation models on non-amenable graphs.
They made a number of conjectures, some of which have since been resolved. 
For a recent review of this field, the reader is referred to the memoir, \cite{Hag11}, by H\"aggstr\"om
on the contributions of Oded Schramm.

The nature of percolation
on a graph $G$ is connected to the algebraic structure of $G$ and its automorphism group.
Of particular interest are Cayley graphs of finitely-generated groups.
The mass-transport principle was introduced to percolation
for graphs with unimodular automorphism group 
by Benjamini, Lyons, Peres, and Schramm, \cite{BLPS};
see \cite[Sect.\ 3]{Hag11}.

\subsection{Disjoint occurrence}

The BK inequality is one of the main tools of percolation theory. 
In its simplest form, it assumes product measure and it applies to increasing events.
The `increasing' assumption was removed by Reimer, \cite{Rei}. The
assumption of independence has proved resilient, but has been weakened by
van den Berg and Jonasson, \cite{vdBJ}, to randomly drawn subsets of given size.
This in turn has provoked a further extension to certain Ising-type models 
by van den Berg
and Gandolfi, \cite{vdBG}.

\subsection{Random-cluster and Potts/Ising models}

Geometrical techniques have long been used in the study of phase transitions.
The random-cluster representation is a prominent example, in which
the two-point correlation functions of the Ising/Potts models are transformed 
into the connection probabilities of a dependent model of percolative type.
Recent books on the random-cluster model include \cite{Grimmett_RCM,Werner_SMF}.

There has been major progress on the theory of random-cluster models, although
many questions remain open. One significant advance concerns the
rigorous theory of the so-called Wulff construction. Wulff, \cite{Wul},
presented a generic formula for the shape of a large `droplet' in
a statistical physical model, as the `crystal' of given volume with smallest aggregate surface 
tension.  Cerf and Pisztora, \cite{CerfPi}, have formulated
and proved the Wulff construction for random-cluster models in general dimension. 
Their solution, together with Cerf's 
solution of the Wulff crystal for percolation, 
is the subject of Rapha\"el Cerf's Saint-Flour
notes, \cite{Cerf-Wu}.

Of the many open problems in this area, we mention one.  
No proof is yet known of the uniqueness
of infinite-volume measures above the critical point in three or more dimensions,
\cite[Thm 5.33]{Grimmett_RCM}.

Just as the two-dimensional Ising model has proved very special, so
are planar random-cluster models with cluster-weighting factor $q=2$.
Smirnov, \cite{Smirnov10}, has shown conformal 
invariance of the critical model on the square lattice, 
and this work has been extended (jointly with Chelkak) to isoradial graphs, 
see \cite{Chelkak-Smirnov3,Chelkak-Smirnov2} and the lecture notes \cite{DS-CI}.
Duminil-Copin, Hongler, and Nolin have extended RSW-type results 
to this model, \cite{DHN}, and this has provided a key technique
for the proof by Lubetzky and Sly, \cite{LS}, of a polynomial
bound for the mixing time of the critical Ising model.  
It is still unclear to what degree techniques for $q=2$ 
may be exploited for general $q \ge 1$.
It has been shown in \cite{BDS} that the so-called `parafermionic observable' of
\cite{Smirnov10} may be used to study the $q>4$ random-cluster model on
the square lattice.

In a further study of the critical Ising model on the re-scaled square lattice $a \ZZ^2$,
Camia, Garban, and Newman, \cite{CGN}, have shown that the magnetization field,
scaled by $a^{15/8}$, converges as $a \to 0$ to a conformally covariant field
which is non-Gaussian. Chelkak, Hongler, and Izuryov, \cite{CHI},
have proved the asymptotic conformal covariance of the $n$-point 
correlation functions of this critical Ising model.
The near-critical FK-Ising and random-cluster models have been considered
by Duminil-Copin, Garban, and Pete, \cite{DGP}.

The exact critical value of the random-cluster model on the square lattice
has long been \lq\lq known but unproven". This has now been rectified by
Beffara and Duminil-Copin, \cite{Beffara_Duminil}, in their proof
that $\pc(q) = \sqrt q/(1+ \sqrt q)$ for $q \ge 1$.
Box-crossings play a central role in their arguments. See also \cite{BDS,G-three}.

\section*{Acknowledgements}
We thank Rapha\"el Cerf,  Michael Damron, Christophe Garban,
and Chuck Newman for their
comments on a draft of this article.
GRG was supported in part by the EPSRC under grant
EP/103372X/1.

\bibliographystyle{amsplain}
\renewcommand{\bibliofont}{\normalsize}
\bibliography{stf2}

\end{document}